\documentclass[12pt]{article}

\newtheorem{proposition}{Proposition}

\newtheorem{theorem}{Theorem}
\newtheorem{corollary}{Corollary}

\newtheorem{remark}{Remark}

\usepackage{amsmath,amssymb,ascmac,color,ulem}

\begin{document}

\title{Basic Invariants of the Complex Reflection Group No.34 Constructed by Conway and Sloane}

\author{Manabu Oura\thanks{Institute of Science and Engineering, Kanazawa University, Japan.\newline
\hspace*{0.5cm}
Email address: oura@se.kanazawa-u.ac.jp} \and Jiro Sekiguchi\thanks{Emeritus Professor, Tokyo University of Agriculture and Technology, Japan.\newline
\hspace*{0.5cm}
Email address: jinomino@yahoo.co.jp
}}



\date{}


\maketitle

\begin{abstract}
This paper studies the basic invariants, constructed by Conway and Sloane, of the complex reflection group
 numbered as 34 in the list of Shephard-Todd \cite{ST}.
\end{abstract}



\section{Introduction}
In this paper, we treat the complex reflection group
of number 34 in the list of Shephard and Todd  \cite{ST} (write this group by ST34 in
this paper)
which has order $6^5\cdot 7!$ and the degrees of its invariants are 6,12,18,24,30,42.
If $x=(x_1,x_2,\cdots,x_6)$ is a linear coordinate of the 6 dimensional representation space of 
ST34, the basic invariants are polynomials of  $x$.
Conway and Sloane \cite{CS} constructed such basic invariants which are denoted by $\mu_j\>(j=6,12,18,24,30,42)$.
It is hard to  write them down  as polynomials of $x$ because of their lengthy.
The main purpose of this paper is to write them in a reasonable size.
We explain our idea to accomplish the purpose briefly.
Let $G(3,3,6)$ be the imprimitive complex reflection group of rank 6.
Since $G(3,3,6)$ is a subgroup of ST34,
each of  $\mu_j\>(j=6,12,18,24,30,42)$ is written as a polynomial of the basic invariants of $G(3,3,6)$,
say $p_3,p_6,p_9,p_{12},p_{15},s_6$.
It is easy to write $p_3,p_6,p_9,p_{12},p_{15},s_6$ down as polynomials of $x$.
Along this idea, we succeeded to write
 $\mu_j\>(j=6,12,18,24,30,42)$ as polynomials of  $p_3,p_6,p_9,p_{12},p_{15},s_6$ in a reasonable size.
Though
it follows from the definition that each $\mu_j$ is invariant by ST34,
we shall give an alternative proof of the invariance by use of $S_6$ invariant polynomials (see Theorem 1).

The discriminant of ST34 is expressed as a polynomial of $\mu_j\>(j=6,12,18,24,30,42)$.
It is known (cf. \cite{BM}) that the discriminant of a complex reflection group  $G$ is
expressed as the determinant of the Saito matrix of $G$.
In ST34 case, Terao and Enta \cite{TE} proposed an algorithm to compute the Saito matrix
for the basic invariants  $\mu_j\>(j=6,12,18,24,30,42)$
and along this line, Bessis and Michel actually computed it explicitly (unpublished).
On the other hand, Kato, Mano and Sekiguchi (\cite{KMS1}, \cite{KMS2}) formulated the notion of the flat structure
which is a generalization of Frobenius manifold structure.
Applying it to a complex reflection group $G$, one can construct the Saito matrix for
the flat coordinate of $G$ when $G$ satisfies  some conditions on invariants.
In ST34 case, the Saito matrix corresponding to a potential vector field of 6 variables
whose weights are the same as the degrees of the basic invariants of ST34 is constructed in \cite{KMS1}
and the determinant of this Saito matrix is expected to be the discriminant
of ST34 if we regard the variables as the set of appropriate basic invariants.
Theorem 2 says that this is actually valid.
Theorem 2 itself is the same as the result in \cite{S2} but its proof is different from
that given in \cite{S2} (see Remark \ref{rem:comment} in the main text).

In the rest of this paper, we shall treat two topics related with the group ST34.
One is the correspondence between the totality of minimal vectors of the Coxeter-Todd lattice and that of 
pseudo-reflections of ST34.
We shall describe a natural $6$-$1$ correspondence between them.
The other is the restriction of the basic invariants $\mu_j\>(j=6,12,18,24,30,42)$
to the representation space of the group ST33 which is the complex reflection group
No.33 in \cite{ST}.
The basic invariants of ST33 were first constructed by H. Burkhardt (cf. \cite{NS}).

This paper is organized as follows.
In section 2, we mention the Conway-Todd lattice, its minimal vectors
and introduce the basic invariants $\mu_j\>(j=6,12,18,24,30,42)$ of ST34 constructed by Conway and Sloane.
In section 3, we first introduce polynomials written as polynomials of
basic invariants $p_3,p_6,p_9,p_{12},p_{15},s_6$ of the complex
reflection group $G(3,3,6)$
which are denoted by $m_j\>(j=6,18,24,30,42)$.
Then in Theorem 1 we shall  show that $m_j$ coincides with $\mu_j$ up to a constant.
In section 4, we mention two results on invariants and a discriminant of ST34
due to Terao, Enta,  Bessis and Michel.
In section 5, we describe the discriminant of ST34 in terms of the flat coordinates of ST34
(cf. \cite{KMS1}, \cite{KMS2}, \cite{S2}).
We include these two sections for the readers who are interested in
the description of  the discriminant
of ST34.
In section 6, we shall treat two topics on ST34.
One is concerned with the correspondence
between the totality of minimal vectors of the Conway-Sloane lattice
and that of hyperplanes fixed by the pseudo-reflections of ST34.
The other is the relationship between the basic invariants of the group ST33 and those of ST34
constructed by Conway and Sloane.
In section 7, we shall mention an application of the result of this paper to deformations
of $E_7$-singularity.
This topic will be discussed elsewhere.

There is a description on a history of basic invariants of the group ST34 in \cite{CS}, \S10.
Concerning the discriminant of ST34,
there are studies by H. Terao and Y. Enta \cite{TE}, P. Orlik and H. Terao \cite{OT}.
Moreover based by the idea of \cite{TE},
D. Bessis and J. Michel computed the Saito matrix of the invariants by Conway and Sloane (unpublished).
Recently Kato, Mano and Sekiguchi \cite{KMS1}, \cite{KMS2} formulated a generalization
of Frobenius manifold structure and  among others they
defined the notion of flat coordinates for well-generated complex reflection groups.
As to the group ST34, the second author (J.S.) constructed the Saito matrix for the flat coordinate.
As a consequence, determinant expression of the discriminant of ST34 for the flat coordinate 
was established (cf. \cite{S2}).

We finally mention that the  softwares Mathematica, Maple and Magma \cite{Bosmaetal} are used to obtain the results in this paper.

\bigskip

\textit{Acknowledgment}.
The second-named author is supported by JSPS KAKENHI Grant Number 17K05269.

\section{The basic invariants  by Conway-Sloane}

By the classification of complex reflection groups of Shephard-Todd \cite{ST},
it is known that there are
 three infinite series of complex reflection groups and thirty-four sporadic ones  numbered as $1,2,\dots,37$.
In this paper  we focus our attention on the  group No.34
which we denote   by ST34.
The basic invariants of ST34 are constructed
in  Theorem 10 of \cite{CS}
and they are denoted by $\mu_j\>(j=6,12,18,24,30,42)$.
We now explain one of their constructions briefly.
Let $\Lambda^{(3)}$ be the Coxeter-Todd lattice 
given by the 3-base (see \cite{CS}). 
Then its automorphism group is nothing else but ST34.
Throughout this paper, we write $\omega = e^{2\pi i/3}$ and  $\theta=\omega-\overline{\omega}=\sqrt{-3}$
without any comment.
The minimal vectors of $\Lambda^{(3)}$ consist of the vectors 

\begin{gather*}
\pm (\omega^a\theta,-\omega^b\theta,0,0,0,0),\\
(\omega^a,\omega^b,\dots, \omega^f), \ a+b+\cdots + f\equiv 0 \pmod{3}
\end{gather*}
where all possible coordinate changes are considered in the first type.
There are $2\cdot 3^2\cdot \binom{6}{2}=270$ minimal vectors of the first type,
whereas $2\cdot 3^5=486$ minimal vectors of the second type.
In total, we have $756$ minimal vectors of $\Lambda^{(3)}$.
We set 
\[
\mu_k=\sum (v_1x_1+\cdots +v_6x_6)^k, \ k=0,1,2\dots,
\]
where $v=(v_1,\dots,v_6)$ runs through $756$ minimal vectors of $\Lambda^{(3)}$.
Theorem 10 in \cite{CS} says that 
the invariant ring of ST34 is generated by 
\[
\mu_6,\mu_{12},\mu_{18},\mu_{24},\mu_{30},\mu_{42}.
\]

We collect the basic properties of ST34.
There are $756$ pseudo-reflections in ST34. 
The  hyperplanes fixed by the pseudo-reflections of ST34 are
\begin{align*}
45 & &  x_i -\omega^a x_j=0,\\
30&&x_1+x_2+x_3+x_4+\omega x_5+\omega^2 x_6=0,\\
20 & & x_1+x_2 + x_3+ \omega (x_4+x_5+x_6)=0,\\
30 & & x_1 +x_2 +\omega (x_3+x_4)+\omega^2 (x_5 +x_6)=0,\\
1 & & x_1 +x_2 +x_3+x_4 +x_5+x_6=0
\end{align*}
where the number on the left side denotes the cardinality of each type.
ST34 can be generated by  $P_1,P_2,P_3,Q_1,R_1,R_2$ given at p.298 in \cite{ST}.
The  hyperplanes fixed by the generators   are
\begin{align*}
P_1 & : 
x_2-x_3=0,\\
P_2&: x_3-x_4=0,\\
P_3&:x_4-x_5=0,\\
Q_1&: x_1-x_2=0,\\
R_1&: x_1-\omega x_2=0,\\
R_2&: x_1+x_2+x_3+x_4+x_5+x_6=0.
\end{align*}
The center $Z_{34}$ of ST34 is of order $6$ generated by 
$-\omega I_6$ where $I_6$ is the identity matrix of degree $6$.

\section{From $G(3,3,6)$ to ST34}
Since the group $G(3,3,6)$ is a subgroup of ST34 and since these two groups are of the same rank,
every invariant by ST34 is a polynomial of the basic invariants of $G(3,3,6)$.

It is known that as the basic invariants of $G(3,3,6)$, we may take
$p_{3j}\>(j=1,2,3,4,5)$ and $s_6$ defined by
$$
p_{3j}=x_1^{3j}+x_2^{3j}+x_3^{3j}+x_4^{3j}+x_5^{3j}+x_6^{3j}\>(j=1,2,3,4,5),
\>s_6=x_1x_2x_3x_4x_5x_6.
$$

Let $m_j\>(j=1,2,3,4,5,7)$ be the polynomials of $p_j\>(j=3,6,9,12,15)$ and $s_6$ defined
by the following identities.
These $m_j$'s are indeed invariants of ST34 as we will show below.

{\footnotesize
$$
\begin{array}{lll}
m_1& =& -5p_3^2 + 6p_6 - 180s_6,\\
m_2 &=&\frac{1}{17}(-10125p_{12} + 1925p_3^4 - 10395p_3^2p_6 + 
     6237p_6^2 + 12375p_3p_9 + 
          51975p_3^2s_6 - 40095p_6s_6 + 935550s_6^2), \\
m_3 &=&\frac{(1}{3644}(-1088916048p_{15}p_3 + 1401597270p_{12}p_3^2 - 
     52286707p_3^6 - 
          3171150p_{12}p_6 + 498415005p_3^4p_6\\
          && - 
     723353895p_3^2p_6^2 + 9550629p_6^3 - 
          961539480p_3^3p_9 + 936512280p_3p_6p_9 - 
     16804260p_9^2 
     + 2817424620p_{12}s_6\\
     && - 
          804053250p_3^4s_6 + 3721617900p_3^2p_6s_6 - 
     1435481190p_6^2s_6 - 
          4300536240p_3p_9s_6
           - 13894040160p_3^2s_6^2 \\
           &&+ 
     8931882960p_6s_6^2 - 
          99578658960s_6^3), \\

m_4& =&\frac{1}{984160}(616868762940p_{12}^2 + 
     13316119811280p_{15}p_3^3 - 16940634571080p_{12}p_3^4 + 
          594592089355p_3^8 \\
          &&- 21567358775952p_{15}p_3p_6 + 
     28072219631940p_{12}p_3^2p_6 - 
          6811339217958p_3^6p_6 - 368171695980p_{12}p_6^2 \\
          &&+ 
     18219805403580p_3^4p_6^2 - 
          14058186924870p_3^2p_6^3 + 98666056665p_6^4 + 
     9579454289856p_{15}p_9\\
     && - 
          13482219282840p_{12}p_3p_9 + 
     11937182241204p_3^5p_9 - 34871337186840p_3^3p_6p_9 
     + 
          24742879985340p_3p_6^2p_9\\
          && + 
     9126295239600p_3^2p_9^2 - 8204834872080p_6p_9^2 + 
          417113860322016p_{15}p_3s_6
           - 
     559068289484760p_{12}p_3^2s_6\\
     &&
      + 24970204832364p_3^6s_6 + 
          17333883551160p_{12}p_6s_6
           - 
     219949713575460p_3^4p_6s_6
      + 
          302052795820380p_3^2p_6^2s_6\\
          && - 
     8679490869060p_6^3s_6
      + 
          401008538842560p_3^3p_9s_6 - 
     376118508300480p_3p_6p_9s_6 + 
          1336769866080p_9^2s_6 \\
          &&
          - 732306533056560p_{12}s_6^2 + 
     311680371562560p_3^4s_6^2 - 
          1202195718884160p_3^2p_6s_6^2
           + 
     386420052498480p_6^2s_6^2 \\
     &&+ 
          1238625892183680p_3p_9s_6^2 + 
     3635549644541760p_3^2s_6^3 
     - 
          1960702790302080p_6s_6^3 + 14488157891797920s_6^4), \\
          \end{array}
          $$
$$
\begin{array}{lll}
m_5&=&\frac{1}{106288160}(-855829791445492800 {p_{12}}^2 {p_3}^2+2089201018632300 {p_{12}}^2 {p_6}
+1165675074483701520 {p_{12}} {p_{15}} {p_3}\\
&&+
138192145506496980 {p_{12}}
   {p_3}^6-877484405515505400 {p_{12}} {p_3}^4 {p_6}
   +1379987318783043000 {p_{12}} {p_3}^3 {p_9}\\
   &&+903656701052469000 {p_{12}} {p_3}^2
   {p_6}^2
   -1000122875204029200 {p_{12}} {p_3} {p_6} {p_9}-
906698664553500 {p_{12}} {p_6}^3\\
&&+580661494716600 {p_{12}} {p_9}^2
-386759386157508288
   {p_{15}}^2-103164357358799784 {p_{15}} {p_3}^5\\
   &&+631029397643881200 {p_{15}} {p_3}^3 {p_6}
   -965383545286463760 {p_{15}} {p_3}^2 {p_9}-
618935794330053240
   {p_{15}} {p_3} {p_6}^2\\
   &&+662798218317617520 {p_{15}} {p_6} {p_9}
   -3721032235241428 {p_3}^{10}+56929735554738525 {p_3}^8 {p_6}\\
   &&-100670803599215340
   {p_3}^7 {p_9}
   -269791048596035610 {p_3}^6 {p_6}^2+754513827309400680 {p_3}^5 {p_6} {p_9}\\
   &&+455456084098977900 {p_3}^4 {p_6}^3
   -540975355278457500
   {p_3}^4 {p_9}^2-1258582004466149700 {p_3}^3 {p_6}^2 {p_9}\\
   &&-237876161937857550 {p_3}^2 {p_6}^4
   +824140299141199800 {p_3}^2 {p_6}
   {p_9}^2+530150266143331200 {p_3} {p_6}^3 {p_9}\\
   &&-157057655983200 {p_3} {p_9}^3
   +168811773953535 {p_6}^5-285007424484525300 {p_6}^2
   {p_9}^2)\\
&&
-\frac{2349}{5314408} (3336381283590 {p_{12}}^2-90991170059715 {p_{12}} {p_3}^4+132473392247970 {p_{12}} {p_3}^2 {p_6}
\\
&&-60473056848750
   {p_{12}} {p_3} {p_9}-3642116633685 {p_{12}} {p_6}^2+69747109280496 {p_{15}} {p_3}^3\\
   &&-96145773383952 {p_{15}} {p_3} {p_6}
   +39062903267592 {p_{15}}
   {p_9}
   +3483778065289 {p_3}^8-38013477459228 {p_3}^6 {p_6}\\
   &&+65401030889973 {p_3}^5 {p_9}+94878525690540 {p_3}^4 {p_6}^2-173998725667650 {p_3}^3
   {p_6} {p_9}\\
   &&-66633212619360 {p_3}^2 {p_6}^3+42238214518500 {p_3}^2 {p_9}^2+110872962731835 {p_3} {p_6}^2 {p_9}+987000196815
   {p_6}^4\\
   &&-32583765485700 {p_6} {p_9}^2)s_6\\
&&
-\frac{221240565}{5314408} (-14731974630 {p_{12}} {p_3}^2+693432090 {p_{12}} {p_6}+10397417904 {p_{15}}
   {p_3}+842268971 {p_3}^6\\
   &&-6659025525 {p_3}^4 {p_6}+11239209720 {p_3}^3 {p_9}+8453844855 {p_3}^2 {p_6}^2-9972788040 {p_3} {p_6}
   {p_9}-348684345 {p_6}^3\\
   &&+86329800 {p_9}^2)s_6^2\\
&&
\frac{98415}{1328602} (7805322022518 {p_{12}}-4961331496841 {p_3}^4+16060185954756 {p_3}^2
   {p_6}-14673183959034 {p_3} {p_9}\\
   &&-4333560899931 {p_6}^2)s_6^3\\
&&
-\frac{231176835}{664301} (9206527958 {p_3}^2-4219165557
   {p_6})s_6^4\\
&&
-\frac{5437195842606312330}{664301}s_6^5,\\
\end{array}
$$

$$
\begin{array}{lll}
m_7&=&
\frac{1}{619872782080}(
-68155702606842557785745 {p_3}^{14}+1432379563699839285857067 {p_6} {p_3}^{12}\\
&&-2557436825642562030969900 {p_9}
   {p_3}^{11}-10698541080341814802394019 {p_6}^2 {p_3}^{10}\\
   &&+3582500028944048652957174 {p_{12}} {p_3}^{10}-2730902888436717092575680
   {p_{15}} {p_3}^9\\
   &&+32708776402929577905944640 {p_6} {p_9} {p_3}^9+34924283387125205511721185 {p_6}^3
   {p_3}^8\\
   &&-25654418121315705916250700 {p_9}^2 {p_3}^8-41680220280677434506011850 {p_{12}} {p_6} {p_3}^8\\
   &&+31369465880845892613346176
   {p_{15}} {p_6} {p_3}^7-123740484999420235384965240 {p_6}^2 {p_9} {p_3}^7\\
   &&+69769318329654728499822240 {p_{12}} {p_9}
   {p_3}^7-52762146031720621233390555 {p_6}^4 {p_3}^6\\
   &&-47257748766395678591117100 {p_{12}}^2 {p_3}^6+126780633297174166548934500
   {p_{12}} {p_6}^2 {p_3}^6\\
   &&+130162353083938619632641600 {p_6} {p_9}^2 {p_3}^6-52307455008394026247756320 {p_{15}} {p_9}
   {p_3}^6\\
   &&-34758789936079286842349760 {p_9}^3 {p_3}^5-92627028581352964521141312 {p_{15}} {p_6}^2
   {p_3}^5\\
   &&+70714813144026703680064512 {p_{12}} {p_{15}} {p_3}^5+186354004242240630334859040 {p_6}^3 {p_9}
   {p_3}^5\\
   &&-268191568191852515860009920 {p_{12}} {p_6} {p_9} {p_3}^5+33299590235273255328597225 {p_6}^5
   {p_3}^4\\
   &&-127770635569767758534767500 {p_{12}} {p_6}^3 {p_3}^4-26421815125729745638761600 {p_{15}}^2
   {p_3}^4\\
   &&-201898474896088367649617400 {p_6}^2 {p_9}^2 {p_3}^4+113298473516080679802360000 {p_{12}} {p_9}^2
   {p_3}^4\\
   &&+122348893127056155271952100 {p_{12}}^2 {p_6} {p_3}^4+190543330015178442011760960 {p_{15}} {p_6} {p_9}
   {p_3}^4\\
   &&+87812086740565495801127040 {p_{15}} {p_6}^3 {p_3}^3+62665614050169799982476800 {p_6} {p_9}^3
   {p_3}^3\\
   &&-74023416302855483255097600 {p_{15}} {p_9}^2 {p_3}^3-167964296041394365510972800 {p_{12}} {p_{15}} {p_6}
   {p_3}^3\\
   &&-106338866778530075623833900 {p_6}^4 {p_9} {p_3}^3-111554371285879313232176400 {p_{12}}^2 {p_9}
   {p_3}^3\\
   &&+262063973509069744833861600 {p_{12}} {p_6}^2 {p_9} {p_3}^3-5126918175165039183578025 {p_6}^6
   {p_3}^2\\
   &&+27931623849212462954949150 {p_{12}} {p_6}^4 {p_3}^2-23152418725131358790400 {p_9}^4 {p_3}^2\\
   &&+29692604870039550095987400
   {p_{12}}^3 {p_3}^2-50201877703761207864827100 {p_{12}}^2 {p_6}^2 {p_3}^2\\
   &&+103065371624630753152920000 {p_6}^3 {p_9}^2
   {p_3}^2-115014969795605791696092000 {p_{12}} {p_6} {p_9}^2 {p_3}^2\\
   &&+56973925496711539138136832 {p_{15}}^2 {p_6}
   {p_3}^2-169930072089259625905241760 {p_{15}} {p_6}^2 {p_9} {p_3}^2\\
   &&+135574470385751288032364160 {p_{12}} {p_{15}} {p_9}
   {p_3}^2-14845481613081683657896320 {p_{15}} {p_6}^4 {p_3}\\
   &&-27945506509393387416667200 {p_6}^2 {p_9}^3
   {p_3}+34586384532087307046400 {p_{12}} {p_9}^3 {p_3}\\
    &&+51824289113657057141222400 {p_{12}} {p_{15}} {p_6}^2
   {p_3}+65877259573354459565260800 {p_{15}} {p_6} {p_9}^2 {p_3}\\
   &&-44231372298721494546912000 {p_{12}}^2 {p_{15}}
   {p_3}+12594866818006221239484000 {p_6}^5 {p_9} {p_3}\\
   &&-44075507458043013535144800 {p_{12}} {p_6}^3 {p_9}
   {p_3}-38836159574544740438311680 {p_{15}}^2 {p_9} {p_3}\\
   &&+37742148440266625813304000 {p_{12}}^2 {p_6} {p_9}
   {p_3}+841205579632152292755 {p_6}^7
   -6381502411588805286450 {p_{12}} {p_6}^5\\
   &&+20837200939432768800000 {p_6}
   {p_9}^4+12591883060088162966100 {p_{12}}^2 {p_6}^3
   -18536934715459277725440 {p_{15}} {p_9}^3\\
&&   +16391367254585859259256064 {p_{12}}
   {p_{15}}^2
   -10659165988259339409769344 {p_{15}}^2 {p_6}^2\\
   &&
   -7683466994926210007484300 {p_6}^4 {p_9}^2
   -8464031877608794242000
   {p_{12}}^2 {p_9}^2\\
&&   +11976082225174674297633600 {p_{12}} {p_6}^2 {p_9}^2+4971887755839338531400 {p_{12}}^3
   {p_6}\\
   &&+18094988745615635242320960 {p_{15}} {p_6}^3 {p_9}-28029510009872734558063680 {p_{12}} {p_{15}} {p_6}
   {p_9})\\
  &&  +
    \frac{2583}{309936391040} (-826650327935043055861 {p_3}^{12}+14425925885422918422282 {p_6} {p_3}^{10}\\
&&    -25452001272005719299540
   {p_9} {p_3}^9-83351888440402255303005 {p_6}^2 {p_3}^8
   +35197343052145422459210 {p_{12}} {p_3}^8\\
   &&-26417038086987608202768
   {p_{15}} {p_3}^7
   +242911215262797035732640 {p_6} {p_9} {p_3}^7+186506835419875414158900 {p_6}^3
   {p_3}^6\\
   &&-183741456439968965505360 {p_9}^2 {p_3}^6-298429989407564825623320 {p_{12}} {p_6} {p_3}^6\\
&&   +220420573635965479389936
   {p_{15}} {p_6} {p_3}^5
   -563000045054149040896440 {p_6}^2 {p_9} {p_3}^5\\
&&   +493130102035689915584400 {p_{12}} {p_9}
   {p_3}^5-161110222055187816091275 {p_6}^4 {p_3}^4\\
&&-328749745772938424739300 {p_{12}}^2 {p_3}^4
+485097352985170111227300 {p_{12}}
   {p_6}^2 {p_3}^4\\
&&   +429094779958543624758000 {p_6} {p_9}^2 {p_3}^4-361895253749323218243360 {p_{15}} {p_9}
   {p_3}^4\\
&&   -20582878620895939008000 {p_9}^3 {p_3}^3
   -347376551409240131304720 {p_{15}} {p_6}^2 {p_3}^3\\
&&   +479349562682466181945440
   {p_{12}} {p_{15}} {p_3}^3
   +428618573732090041704000 {p_6}^3 {p_9} {p_3}^3\\
&&   -673892039580382647561600 {p_{12}} {p_6}
   {p_9} {p_3}^3+36297779413268324630250 {p_6}^5 {p_3}^2\\
 &&  -143166572523976223269800 {p_{12}} {p_6}^3
   {p_3}^2-173568570666357873223296 {p_{15}}^2 {p_3}^2\\
 &&  -293344069145812931192400 {p_6}^2 {p_9}^2 {p_3}^2+43066790405325599148000
   {p_{12}} {p_9}^2 {p_3}^2\\
 &&  +141143057326507215360600 {p_{12}}^2 {p_6} {p_3}^2+469676819583053026119360 {p_{15}} {p_6}
   {p_9} {p_3}^2\\
   &&+90703842698655568839600 {p_{15}} {p_6}^3 {p_3}+19163765066700626496000 {p_6} {p_9}^3
   {p_3}\\
 &&  -22962898216765380119040 {p_{15}} {p_9}^2 {p_3}
   -178484620941100228128480 {p_{12}} {p_{15}} {p_6}
   {p_3}\\
 \end{array}
   $$
   $$
   \begin{array}{lll}
 &&  -81029269710700131377700 {p_6}^4 {p_9} {p_3}-19161435544300302037200 {p_{12}}^2 {p_9} {p_3}\\
&&   +169202690869163631008400
   {p_{12}} {p_6}^2 {p_9} {p_3}-44454150584215675875 {p_6}^6
   +236755439504502230250 {p_{12}} {p_6}^4\\
   &&-95534751065126400
   {p_9}^4+236549851238374857000 {p_{12}}^3
   -413969210349675356700 {p_{12}}^2 {p_6}^2\\
   &&+43860401994452195041200 {p_6}^3
   {p_9}^2
   -11869779935007266508000 {p_{12}} {p_6} {p_9}^2+53589954112307329756032 {p_{15}}^2 {p_6}\\
   &&-97278586484753853185760
   {p_{15}} {p_6}^2 {p_9}+14219410871098234843200 {p_{12}} {p_{15}} {p_9})s_6\\
  && -\frac{637362999}{15496819552}
  (4857452670185866 {p_3}^{10}-65282378295761109 {p_6} {p_3}^8+108907344194688000 {p_9} {p_3}^7\\
 && +260295982554842820 {p_6}^2
   {p_3}^6-144443551320429948 {p_{12}} {p_3}^6+105268487601706560 {p_{15}} {p_3}^5\\
   &&-644181269743573344 {p_6} {p_9}
   {p_3}^5-367603906096725030 {p_6}^3 {p_3}^4+384498967408122240 {p_9}^2 {p_3}^4\\
  && +678555020839017240 {p_{12}} {p_6}
   {p_3}^4-463955787849692736 {p_{15}} {p_6} {p_3}^3
   +895430833203929760 {p_6}^2 {p_9} {p_3}^3\\
   &&-880017911185364640
   {p_{12}} {p_9} {p_3}^3+153022160508667650 {p_6}^4 {p_3}^2
   +472031454180039120 {p_{12}}^2 {p_3}^2\\
   &&-541859921711888940
   {p_{12}} {p_6}^2 {p_3}^2-505515465619816320 {p_6} {p_9}^2 {p_3}^2
   +568858413963720960 {p_{15}} {p_9}
   {p_3}^2\\
   &&+1078387051311360 {p_9}^3 {p_3}+333737005975634880 {p_{15}} {p_6}^2 {p_3}
   -573606450773307072 {p_{12}} {p_{15}}
   {p_3}\\
   &&-298809060268406400 {p_6}^3 {p_9} {p_3}+519137518908140640 {p_{12}} {p_6} {p_9} {p_3}
   -740481396190785
   {p_6}^5\\
   &&+2897609287598280 {p_{12}} {p_6}^3
   +163131844841144064 {p_{15}}^2
   +133116924334917600 {p_6}^2 {p_9}^2\\
   &&-769650476421600
   {p_{12}} {p_9}^2-2833333167373020 {p_{12}}^2 {p_6}
   -295206239616621696 {p_{15}} {p_6} {p_9})s_6^2\\
&& -\frac{137781}{7748409776}
   (212849412401280527467 {p_3}^8-2056737692333330079912 {p_6} {p_3}^6
   +3301675794233514997830 {p_9}
   {p_3}^5\\
   &&+4487246757920672916210 {p_6}^2 {p_3}^4-4310703736169475684240 {p_{12}} {p_3}^4
   +3103075441442929625568 {p_{15}}
   {p_3}^3\\
   &&-7542458411263968569940 {p_6} {p_9} {p_3}^3
   -2624550943988148988320 {p_6}^3 {p_3}^2+1648635226201837343160
   {p_9}^2 {p_3}^2\\
   &&+5112301708617615116280 {p_{12}} {p_6} {p_3}^2-3259289326066664689728 {p_{15}} {p_6}
   {p_3}
   +3870577632979310328810 {p_6}^2 {p_9} {p_3}\\
   &&-2142794233908947745900 {p_{12}} {p_9} {p_3}+66893663401936554375
   {p_6}^4
   +193991306712427220580 {p_{12}}^2\\
   &&-230654127137064281640 {p_{12}} {p_6}^2-903663040115630763720 {p_6}
   {p_9}^2
   +1073610002137218186960 {p_{15}} {p_9})s_6^3\\
   &&  
   -\frac{457570701 }{3874204888} (393541072937010113 {p_3}^6-2483580961309967805
   {p_6} {p_3}^4+3573657602659051360 {p_9} {p_3}^3\\
  && +2671680990475403505 {p_6}^2 {p_3}^2
   -4136783742478404450 {p_{12}}
   {p_3}^2+2632941690092063712 {p_{15}} {p_3}\\
   &&-2833712226159225840 {p_6} {p_9} {p_3}-136991407212027225
   {p_6}^3+52735937047941600 {p_9}^2
   +266860623034204470 {p_{12}} {p_6})s_6^4\\
  && +  
   \frac{125452159093170 }{484275611}(-1416246125627
   {p_3}^4+3374374908123 {p_6} {p_3}^2-2450724434568 {p_9} {p_3}\\
   &&-673355559033 {p_6}^2+1069099330203
   {p_{12}})s_6^5\\
 &&  
     -\frac{9039811410 }{484275611}(105868047331971014 {p_3}^2-36630825385046211
   {p_6})s_6^6\\
   &&
   -\frac{1266409465981399253335790610}{484275611}s_6^7.\\
   \end{array}
   $$
}

The following theorem is the main result of this paper.
\begin{theorem}
 The polynomials $m_j\>(j=1,2,3,4,5,7)$ are invariants of ST34.
\end{theorem}

\textit{Proof}.
The group ST34 can be generated by the pseudo-reflections $P_1,P_2,P_3,Q_1,R_1,R_2$.
It is clear that $p_{3j}\>(j=1,2,3,4,5)$ and $s_6$ are invariant under the action 
of $P_1,P_2,P_3,Q_1,R_1$ and so are $m_j\>(j=1,2,3,4,5,7)$.

As a consequence, it is sufficient
to show that each of $m_j$ is invariant by the action of $R_2$.
For this purpose, we  introduce the polynomials
$$
q_j=\sum_{k=1}^6x_k^j\quad
(j=1,2,3,4,5).
$$
Then 
$p_{3j}\>(j=1,2,3,4,5)$  are polynomials of $q_j\>(j=1,2,3,4,5)$ and $s_6$.
Indeed, we have
{\footnotesize
$$
\begin{array}{lll}
p_3& =&  q_3, \\
  p_6 &=& \frac{1}{120}  ( q_1^6 - 15   q_1^4   q_2 + 45   q_1^2   q_2^2 - 15   q_2^3 + 
     40   q_1^3   q_3 - 120   q_1   q_2   q_3 + 40   q_3^2 - 90   q_1^2   q_4 + 
     90   q_2   q_4 + 
          144   q_1   q_5 - 720    s_6), \\
  p_9& =& \frac{1}{720}  ( q_1^9 - 9   q_1^7   q_2 - 27   q_1^5   q_2^2 + 
     135   q_1^3   q_2^3 + 33   q_1^6   q_3 + 45   q_1^4   q_2   q_3 - 
     135   q_1^2   q_2^2   q_3 - 
          135   q_2^3   q_3 + 120   q_1^3   q_3^2 - 360   q_1   q_2   q_3^2\\
          && + 
     80   q_3^3 - 54   q_1^5   q_4 - 270   q_1^3   q_2   q_4 - 
     270   q_1^2   q_3   q_4 + 270   q_2   q_3   q_4 + 54   q_1^4   q_5 + 
          324   q_1^2   q_2   q_5 + 162   q_2^2   q_5 + 432   q_1   q_3   q_5 \\
          &&+ 
     324   q_4   q_5 - 1080   q_1^3    s_6 - 3240   q_1   q_2    s_6 - 
     2160   q_3    s_6), \\
    p_{12} &= &\frac{1}{43200}  ( q_1^{12} - 12   q_1^{10 }  q_2 + 90   q_1^8   q_2^2 - 
     900   q_1^6   q_2^3 + 2025   q_1^4   q_2^4 + 80   q_1^9   q_3 - 
     720   q_1^7   q_2   q_3 + 3600   q_1^5   q_2^2   q_3 - 
          3600   q_1^3   q_2^3   q_3 \\
          &&+ 1680   q_1^6   q_3^2 - 
     3600   q_1^4   q_2   q_3^2 + 3600   q_1^2   q_2^2   q_3^2 - 
     3600   q_2^3   q_3^2 + 3200   q_1^3   q_3^3 - 9600   q_1   q_2   q_3^3 + 
          1600   q_3^4 - 45   q_1^8   q_4 + 180   q_1^6   q_2   q_4\\
          && - 
     6750   q_1^4   q_2^2   q_4 + 8100   q_1^2   q_2^3   q_4 - 2025   q_2^4   q_4 - 
     2880   q_1^5   q_3   q_4 + 7200   q_1^3   q_2   q_3   q_4 - 
          21600   q_1   q_2^2   q_3   q_4 - 7200   q_1^2   q_3^2   q_4\\
          && + 
     7200   q_2   q_3^2   q_4 - 16200   q_1^2   q_2   q_4^2 + 8100   q_2^2   q_4^2 + 
     2700   q_4^3 + 288   q_1^7   q_5 - 
          3456   q_1^5   q_2   q_5 + 12960   q_1^3   q_2^2   q_5
           + 
     11520   q_1^4   q_3   q_5 \\
     &&- 17280   q_1^2   q_2   q_3   q_5 + 
     11520   q_1   q_3^2   q_5 - 17280   q_1^3   q_4   q_5 + 
          25920   q_1   q_2   q_4   q_5 + 17280   q_3   q_4   q_5 + 
     20736   q_1^2   q_5^2 + 10368   q_2   q_5^2\\
     && - 1440   q_1^6    s_6 - 
     64800   q_1^2   q_2^2    s_6 - 57600   q_1^3   q_3    s_6 - 
          57600   q_3^2    s_6 - 129600   q_2   q_4    s_6 - 
     207360   q_1   q_5    s_6 + 259200    s_6^2), \\
    p_{15} &=&
     \frac{1}{1036800} (13   q_1^{15} - 330   q_1^{13 }  q_2 + 
     2565   q_1^{11}   q_2^2 - 4500   q_1^9   q_2^3 - 15525   q_1^7   q_2^4 + 
     44550   q_1^5   q_2^5 - 10125   q_1^3   q_2^6 + 
          895   q_1^{12}   q_3 \\
          &&- 12090   q_1^{10}   q_2   q_3 + 
     22275   q_1^8   q_2^2   q_3 + 78300   q_1^6   q_2^3   q_3 - 
     165375   q_1^4   q_2^4   q_3 - 20250   q_1^2   q_2^5   q_3 + 
     10125   q_2^6   q_3 + 
          16000   q_1^9   q_3^2\\
          && - 36000   q_1^7   q_2   q_3^2 - 
     64800   q_1^5   q_2^2   q_3^2 + 108000   q_1^3   q_2^3   q_3^2 + 
     108000   q_1   q_2^4   q_3^2 + 24800   q_1^6   q_3^3 - 
          84000   q_1^4   q_2   q_3^3 + 108000   q_1^2   q_2^2   q_3^3 \\
          &&- 
     36000   q_2^3   q_3^3 + 32000   q_1^3   q_3^4 - 96000   q_1   q_2   q_3^4 + 
     12800   q_3^5 - 1710   q_1^{11}   q_4 + 
          19800   q_1^9   q_2   q_4 - 2700   q_1^7   q_2^2   q_4 \\
          &&-
     178200   q_1^5   q_2^3   q_4 
     + 101250   q_1^3   q_2^4   q_4 - 
     58050   q_1^8   q_3   q_4 - 16200   q_1^6   q_2   q_3   q_4 + 
          337500   q_1^4   q_2^2   q_3   q_4 + 81000   q_1^2   q_2^3   q_3   q_4\\
          && - 
     101250   q_2^4   q_3   q_4 - 43200   q_1^5   q_3^2   q_4 + 
     216000   q_1^3   q_2   q_3^2   q_4 - 
          432000   q_1   q_2^2   q_3^2   q_4 - 72000   q_1^2   q_3^3   q_4 + 
     72000   q_2   q_3^3   q_4\\
     && + 56700   q_1^7   q_4^2 + 
     113400   q_1^5   q_2   q_4^2 - 121500   q_1^3   q_2^2   q_4^2 + 
          40500   q_1^4   q_3   q_4^2 - 405000   q_1^2   q_2   q_3   q_4^2 + 
     121500   q_2^2   q_3   q_4^2 - 81000   q_1^3   q_4^3 \\
     &&+ 81000   q_3   q_4^3 + 
     2196   q_1^{10}   q_5 - 22140   q_1^8   q_2   q_5 - 
          33480   q_1^6   q_2^2   q_5 + 243000   q_1^4   q_2^3   q_5 + 
     72900   q_1^2   q_2^4   q_5 - 24300   q_2^5   q_5 \\
     &&+ 77760   q_1^7   q_3   q_5 + 
     60480   q_1^5   q_2   q_3   q_5 - 
          302400   q_1^3   q_2^2   q_3   q_5 - 388800   q_1   q_2^3   q_3   q_5 + 
     216000   q_1^4   q_3^2   q_5 
     - 432000   q_1^2   q_2   q_3^2   q_5\\
     && - 
     43200   q_2^2   q_3^2   q_5 + 115200   q_1   q_3^3   q_5 - 
          136080   q_1^6   q_4   q_5 - 550800   q_1^4   q_2   q_4   q_5 + 
     97200   q_1^2   q_2^2   q_4   q_5 + 97200   q_2^3   q_4   q_5 \\
     &&- 
     259200   q_1^3   q_3   q_4   q_5 + 259200   q_1   q_2   q_3   q_4   q_5 + 
          259200   q_3^2   q_4   q_5 - 97200   q_1^2   q_4^2   q_5 + 
     291600   q_2   q_4^2   q_5 + 82944   q_1^5   q_5^2\\
     && + 
     414720   q_1^3   q_2   q_5^2 + 311040   q_1   q_2^2   q_5^2 + 
          414720   q_1^2   q_3   q_5^2 + 207360   q_2   q_3   q_5^2 + 
     311040   q_1   q_4   q_5^2 + 41472   q_5^3 - 23400   q_1^9    s_6 \\
     &&+ 
     280800   q_1^7   q_2    s_6 - 226800   q_1^5   q_2^2    s_6 - 
          1620000   q_1^3   q_2^3    s_6 + 243000   q_1   q_2^4    s_6 - 
     820800   q_1^6   q_3    s_6 + 216000   q_1^4   q_2   q_3    s_6 \\
     &&+ 
     3240000   q_1^2   q_2^2   q_3    s_6 + 
          648000   q_2^3   q_3    s_6 - 864000   q_1^3   q_3^2    s_6 + 
     864000   q_1   q_2   q_3^2    s_6 - 576000   q_3^3    s_6
      + 
     1555200   q_1^5   q_4    s_6 \\
     &&+ 2592000   q_1^3   q_2   q_4    s_6 - 
          1944000   q_1   q_2^2   q_4    s_6 - 2592000   q_2   q_3   q_4    s_6 - 
     972000   q_1   q_4^2    s_6 - 2073600   q_1^4   q_5    s_6 \\
     &&- 
     6220800   q_1^2   q_2   q_5    s_6 - 
          1555200   q_2^2   q_5    s_6 - 4147200   q_1   q_3   q_5    s_6 - 
     1555200   q_4   q_5    s_6 + 10368000   q_1^3    s_6^2 + 
     15552000   q_1   q_2    s_6^2\\
     && + 5184000   q_3    s_6^2).
     \\
     \end{array}
     $$
}
The reflection $R_2$ acts on $(x_1,x_2,x_3,x_4,x_5,x_6)$ by
$$
R_2:x_j\mapsto x_j-\frac{1}{3}\sum_{j=k}^6x_k\quad(j=1,2,\cdots,6).
$$
If  $f=f(x_1,\cdots,x_6)$ is a polynomial of $(x_1,\cdots,x_6)$,
we write $f\circ R_2$ for $f(R_2(x_1,\cdots,x_6))$.
Then by direct computation, we have
{\footnotesize
$$
\begin{array}{lll}
p_1\circ R_2 &=&\frac{1}{9} (q_1^3 - 9 q_1 q_2 + 9 q_3), \\
 p_2\circ R_2 &=&\frac{1}{9720} (-79 q_1^6 + 585 q_1^4 q_2 + 3645 q_1^2 q_2^2 - 1215 q_2^3 - 
     3960 q_1^3 q_3 - 9720 q_1 q_2 q_3 + 3240 q_3^2 + 
          8910 q_1^2 q_4 + 7290 q_2 q_4 \\
          &&- 7776 q_1 q_5 - 58320 s_6),\\
    p_3\circ R_2 &=& \frac{1}{524880} (-1384 q_1^9 + 22707 q_1^7 q_2 - 
     85293 q_1^5 q_2^2 - 91125 q_1^3 q_2^3 + 98415 q_1 q_2^4 - 
     57375 q_1^6 q_3 + 441045 q_1^4 q_2 q_3\\
     && + 
          338985 q_1^2 q_2^2 q_3 - 98415 q_2^3 q_3 - 
     456840 q_1^3 q_3^2 - 262440 q_1 q_2 q_3^2 + 58320 q_3^3 + 
     82134 q_1^5 q_4 - 634230 q_1^3 q_2 q_4 \\
     &&- 
          393660 q_1 q_2^2 q_4 + 1027890 q_1^2 q_3 q_4 + 
     196830 q_2 q_3 q_4 - 393660 q_1 q_4^2 - 53946 q_1^4 q_5 + 
     446148 q_1^2 q_2 q_5 + 118098 q_2^2 q_5\\
     && - 
          524880 q_1 q_3 q_5 + 236196 q_4 q_5 + 612360 q_1^3 s_6 + 
     3936600 q_1 q_2 s_6 - 1574640 q_3 s_6), \\
   p_4\circ R_2 &=&\frac{1}{283435200} (-61441 q_1^{12} + 1761678 q_1^{10} q_2 - 
     17211690 q_1^8 q_2^2 + 57445200 q_1^6 q_2^3 + 
     2460375 q_1^4 q_2^4 - 32476950 q_1^2 q_2^5\\
     && - 
          2696580 q_1^9 q_3 + 59428080 q_1^7 q_2 q_3 - 
     327612600 q_1^5 q_2^2 q_3 - 23619600 q_1^3 q_2^3 q_3 + 
     64953900 q_1 q_2^4 q_3 - 34739280 q_1^6 q_3^2\\
     && + 
          553748400 q_1^4 q_2 q_3^2 + 81356400 q_1^2 q_2^2 q_3^2 - 
     23619600 q_2^3 q_3^2 - 235612800 q_1^3 q_3^3 - 
     62985600 q_1 q_2 q_3^3 + 10497600 q_3^4 \\
     &&+ 
          1843155 q_1^8 q_4 - 58480380 q_1^6 q_2 q_4 + 
     388739250 q_1^4 q_2^2 q_4 + 53144100 q_1^2 q_2^3 q_4 - 
     13286025 q_2^4 q_4 + 19595520 q_1^5 q_3 q_4\\
     && - 
          1165233600 q_1^3 q_2 q_3 q_4 - 
     141717600 q_1 q_2^2 q_3 q_4 + 530128800 q_1^2 q_3^2 q_4 + 
     47239200 q_2 q_3^2 q_4 + 43302600 q_1^4 q_4^2\\
     && + 
          543250800 q_1^2 q_2 q_4^2 + 53144100 q_2^2 q_4^2 - 
     259815600 q_1 q_3 q_4^2 + 17714700 q_4^3 + 93312 q_1^7 q_5 + 
     23514624 q_1^5 q_2 q_5\\
     && - 
          226748160 q_1^3 q_2^2 q_5 + 29393280 q_1^4 q_3 q_5 + 
     579467520 q_1^2 q_2 q_3 q_5 - 201553920 q_1 q_3^2 q_5 - 
     113374080 q_1^3 q_4 q_5 \\
     &&- 
          453496320 q_1 q_2 q_4 q_5 + 113374080 q_3 q_4 q_5 + 
     52907904 q_1^2 q_5^2 + 68024448 q_2 q_5^2 + 
     75232800 q_1^6 s_6 \\
     &&- 692841600 q_1^4 q_2 s_6 - 
          1464415200 q_1^2 q_2^2 s_6 + 1700611200 q_1^3 q_3 s_6 + 
     2078524800 q_1 q_2 q_3 s_6 - 377913600 q_3^2 s_6\\
     && - 
     2078524800 q_1^2 q_4 s_6 - 
          850305600 q_2 q_4 s_6 + 1133740800 q_1 q_5 s_6 + 
     1700611200 s_6^2), \\
   p_5\circ R_2 &=& \frac{1}{61222003200} (190052 q_1^{15} + 7543305 q_1^{13} q_2 - 
     314760330 q_1^{11} q_2^2 + 3315388725 q_1^9 q_2^3 - 
     11699575200 q_1^7 q_2^4 \\
     &&+ 
          4687014375 q_1^5 q_2^5 + 3587226750 q_1^3 q_2^6 - 
     597871125 q_1 q_2^7 + 41009175 q_1^{12} q_3 + 
     135576990 q_1^{10} q_2 q_3\\
     && - 12521395125 q_1^8 q_2^2 q_3 + 
          74175385500 q_1^6 q_2^3 q_3 - 20512146375 q_1^4 q_2^4 q_3 - 
     12356003250 q_1^2 q_2^5 q_3 + 597871125 q_2^6 q_3 \\
     &&+ 
     1263778200 q_1^9 q_3^2 + 
          7463793600 q_1^7 q_2 q_3^2 - 
     163560481200 q_1^5 q_2^2 q_3^2 + 17950896000 q_1^3 q_2^3 q_3^2 + 
     10097379000 q_1 q_2^4 q_3^2\\
     && + 6796029600 q_1^6 q_3^3 + 
          138332124000 q_1^4 q_2 q_3^3 + 
     3070548000 q_1^2 q_2^2 q_3^3 - 2125764000 q_2^3 q_3^3 - 
     36321696000 q_1^3 q_3^4 \\
     &&- 5668704000 q_1 q_2 q_3^4 + 
     755827200 q_3^5 - 
          124608780 q_1^{11} q_4 + 1213748550 q_1^9 q_2 q_4 + 
     7755102000 q_1^7 q_2^2 q_4\\
     && - 79344141300 q_1^5 q_2^3 q_4 + 
     12488863500 q_1^3 q_2^4 q_4 + 
          8370195750 q_1 q_2^5 q_4 - 6226279650 q_1^8 q_3 q_4 + 
     17065161000 q_1^6 q_2 q_3 q_4\\
     && + 
     320016055500 q_1^4 q_2^2 q_3 q_4 - 
     2657205000 q_1^2 q_2^3 q_3 q_4 - 
          5978711250 q_2^4 q_3 q_4 - 55642528800 q_1^5 q_3^2 q_4 \\
          &&- 
     327840048000 q_1^3 q_2 q_3^2 q_4 - 
     10628820000 q_1 q_2^2 q_3^2 q_4 + 81723816000 q_1^2 q_3^3 q_4 + 
          4251528000 q_2 q_3^3 q_4\\
          && + 5558479200 q_1^7 q_4^2 - 
     30711384900 q_1^5 q_2 q_4^2 - 150397803000 q_1^3 q_2^2 q_4^2 - 
     16740391500 q_1 q_2^3 q_4^2 + 
          90020200500 q_1^4 q_3 q_4^2\\
          && + 
     214170723000 q_1^2 q_2 q_3 q_4^2 + 7174453500 q_2^2 q_3 q_4^2 - 
     44641044000 q_1 q_3^2 q_4^2 - 35783694000 q_1^3 q_4^3\\
     && - 
          33480783000 q_1 q_2 q_4^3 + 4782969000 q_3 q_4^3 + 
     122974524 q_1^{10} q_5 - 1648516860 q_1^8 q_2 q_5 + 
     569232360 q_1^6 q_2^2 q_5\\
     && + 41133533400 q_1^4 q_2^3 q_5 - 
          11319693300 q_1^2 q_2^4 q_5 - 1434890700 q_2^5 q_5 + 
     5427259200 q_1^7 q_3 q_5 - 30377954880 q_1^5 q_2 q_3 q_5\\
     && - 
     148803480000 q_1^3 q_2^2 q_3 q_5 + 
          12754584000 q_1 q_2^3 q_3 q_5 + 
     45381124800 q_1^4 q_3^2 q_5 + 125278358400 q_1^2 q_2 q_3^2 q_5\\
     && - 
     2550916800 q_2^2 q_3^2 q_5 - 24942297600 q_1 q_3^3 q_5 - 
          8299927440 q_1^6 q_4 q_5 + 58741945200 q_1^4 q_2 q_4 q_5 \\
          &&+ 
     130734486000 q_1^2 q_2^2 q_4 q_5 + 5739562800 q_2^3 q_4 q_5 - 
     118475913600 q_1^3 q_3 q_4 q_5 - 
          127545840000 q_1 q_2 q_3 q_4 q_5\\
          && + 
     15305500800 q_3^2 q_4 q_5 + 56757898800 q_1^2 q_4^2 q_5 + 
     17218688400 q_2 q_4^2 q_5 + 2781444096 q_1^5 q_5^2 - 
          23128312320 q_1^3 q_2 q_5^2 \\
          &&- 24488801280 q_1 q_2^2 q_5^2 + 
     34012224000 q_1^2 q_3 q_5^2 + 12244400640 q_2 q_3 q_5^2 - 
     24488801280 q_1 q_4 q_5^2 + 
          2448880128 q_5^3\\
          && + 1452897000 q_1^9 s_6 - 
     35980524000 q_1^7 q_2 s_6 + 213119650800 q_1^5 q_2^2 s_6 + 
     172186884000 q_1^3 q_2^3 s_6 \\
     &&- 
          52612659000 q_1 q_2^4 s_6 + 45885009600 q_1^6 q_3 s_6 - 
     807317928000 q_1^4 q_2 q_3 s_6 - 
     522937944000 q_1^2 q_2^2 q_3 s_6\\
     && + 
          38263752000 q_2^3 q_3 s_6 + 
     530968608000 q_1^3 q_3^2 s_6 + 
     289103904000 q_1 q_2 q_3^2 s_6 - 34012224000 q_3^3 s_6 \\
     &&- 
     20596291200 q_1^5 q_4 s_6 + 
          926833104000 q_1^3 q_2 q_4 s_6 + 
     420901272000 q_1 q_2^2 q_4 s_6 - 
     833299488000 q_1^2 q_3 q_4 s_6 \\
     &&- 
     153055008000 q_2 q_3 q_4 s_6 + 
          210450636000 q_1 q_4^2 s_6 - 
     11337408000 q_1^4 q_5 s_6 - 510183360000 q_1^2 q_2 q_5 s_6 \\
     &&- 
     91833004800 q_2^2 q_5 s_6 + 326517350400 q_1 q_3 q_5 s_6 - 
          91833004800 q_4 q_5 s_6 - 578207808000 q_1^3 s_6^2 \\
          &&- 
     1224440064000 q_1 q_2 s_6^2 + 306110016000 q_3 s_6^2), \\
   s_6\circ R_2 &=& (-13 q_1^6 + 180 q_1^4 q_2 - 405 q_1^2 q_2^2 - 
     450 q_1^3 q_3 + 810 q_1 q_2 q_3 + 810 q_1^2 q_4 - 972 q_1 q_5 + 
     14580 s_6)/14580.\\
     \end{array}
     $$
}
We now take one of 
$m_1,m_2,m_3,m_4,m_5,m_7$ and write it, say $m_j$.
It is possible to show the expressions of  both $m_j$ and $m_j\circ R_2$
as polynomials of $q_1,q_2,q_3,q_4,q_5$ and $s_6$.
Comparing these two expressions, we find that $m_j\circ R_2=m_j$.
This completes the proof of the theorem.

\begin{remark}
The normalization of the choice of $m_j$'s is that as a polynomial of $x$, each $m_j$ contains
the term $x_1^{6j}$ with coefficient $1$. This implies that as a polynomial of
$p_3,p_6,p_9,p_{12},p_{15},s_6$, the value of $m_j$ at $p_3=p_6=p_9=p_{12}=p_{15}=1,s_6=0$ is $1$.
\end{remark}

\begin{corollary}
The relations between $\mu_6,\mu_{12},\cdots$ and $m_1,m_2,\cdots$
are given by
$$
\begin{array}{lll}
\mu_6&=&  -1944m_1,\\
\mu_{12}&=&  66096 m_2,\\
\mu_{18}&=&  -1770984 m_3,\\
\mu_{24}&=&  47830176m_4,\\
\mu_{30}&=&  -1291401144 m_5,\\
\mu_{42}&=&  -941431787784 m_7.\\
\end{array}
$$
\end{corollary}

\textit{Proof}.
The polynomials $m_j\>(j=1,2,3,4,5,7)$ are obtained by the same idea 
 as the construction of $\mu_{6j}\>(j=1,2,3,4,5,7)$ given in \cite{CS}.
Actually $m_j$ coincides with $\mu_{6j}$ up to a constant factor.
The constants 
are specified by
evaluating polynomials at special values and the result follows.

\section{Two results on invariants and a discriminant of ST34}

We give two results which play important roles in our consideration.

\subsection{Invariants by Terao and Enta}

We recall the invariants $f_1,f_2,f_3,f_4,f_5,f_6$ given in Terao and Enta \cite{TE}
(see also \cite{OT}).
By direct computation we have
$$
\begin{array}{lll}
f_1&=&\frac{1}{1944}\mu_6,\\
f_2&=&\frac{1}{3888}\mu_{12},\\
f_3&=&\frac{1}{1944}\mu_{18},\\
f_4&=&\frac{1937160963}{1797549546875}m_1^4 -\frac{31670896436}{
    19773045015625}m_1^2m_2 +\frac{233872961}{
    754856437500}m_2^2 +\frac{63038467}{258156875724} m_1m_3- 
      \frac{6151}{205320951}m_{4},\\
      f_5&=&\frac{1}{1944}\mu_{30},\\
      f_6&=&\frac{1}{1944}\mu_{42},\\
      \end{array}
      $$
where $f_4$ is a multiple of the Hessian of $f_1$.

       \subsection{Discriminant by Bessis and Michel}
       
Concerning the discriminant of ST34,
there are studies by H. Terao and Y. Enta (\cite{TE},  \cite{OT}).
D. Bessis and J. Michel \cite{BM} constructed Saito matrices of some complex reflection groups.
They commented that the case ST34 is not treated in \cite{BM} since the Saito matrix of ST34 is too long
to include it in the paper.
Moreover based on the idea of \cite{TE},
they constructed the Saito matrix of the invariants by Conway and Sloane
explicitly (unpublished).

Prof. J. Michel kindly sent the data of the Saito matrix of ST34 for the basic invariants
$\{f_1,f_2,\cdots,f_6\}$ on the request of the second author (J. S.).
We call this matrix $M_{34}$ whose determinant is the discriminant of ST34
up to a non-zero constant factor.
Each matrix  entry of $M_{34}$ is a polynomial of $x,y,z,t,u,v$ which are the same as $f_1,\dots,f_6$ of Terao and Enta
in this ordering.

\section{Construction of the Saito matrix in terms of a potential vector field}

The Saito matrix of a reflection group depends on the choice of basic invariants.
There is an alternative way to construct the Saito matrix of ST34.
We begin this section with explaining the construction of the Saito matrix for the flat coordinate.
For the details, see  \cite{S2}, \cite{KMS1}, \cite{S1}.

The Saito matrix of ST34 can be expressed by a potential vector field
given in \cite{KMS1}.
To show the result, we define the following polynomials $h_j\>(j=1,2,3,4,5,6)$ of $u_1,u_2,u_3,u_4,u_5,u_6$.

{\footnotesize
$$
\begin{array}{ll}
h_1&\hspace{-2.5mm}=(20u_1^4u_2^2 + 120u_1^2u_2^3 - 60u_2^4 - 12u_1^5u_3 + 
    60u_1^3u_2u_3
 - 180u_1u_2^2u_3 + 135u_1^2u_3^2 + 135u_2u_3^2 + 
    120u_1^4u_4 \\
&\hspace{2mm}- 
        180u_1^2u_2u_4 + 540u_2^2u_4 + 270u_1u_3u_4 + 405u_4^2 + 
    180u_1^3u_5 + 540u_1u_2u_5 + 405u_3u_5 + 405u_1u_6)/405, \\
h_2&\hspace{-2.5mm}= (64u_1^9 - 288u_1^7u_2 - 1728u_1^5u_2^2 + 
    1728u_1^3u_2^3 - 2592u_1u_2^4 + 432u_1^6u_3
 + 3888u_1^4u_2u_3 + 
    3888u_1^2u_2^2u_3 \\
&\hspace{2mm}+ 
        1944u_2^3u_3 + 972u_1^3u_3^2 + 729u_3^3 - 1296u_1^5u_4 + 
    7776u_1^3u_2u_4 + 8748u_1^2u_3u_4 - 2916u_2u_3u_4 - 
    5832u_1u_4^2\\
&\hspace{2mm} + 3888u_1^4u_5 - 
        2916u_1^2u_2u_5 - 2916u_2^2u_5 + 4374u_1u_3u_5 + 4374u_4u_5 + 
    4374u_2u_6)/4374, \\
 h_3&\hspace{-2.5mm}= (64u_1^{10} + 3456u_1^8u_2 + 6624u_1^6u_2^2 + 
    2160u_1^4u_2^3 + 3888u_2^5 + 7776u_1^7u_3 + 19440u_1^5u_2u_3 + 
    32400u_1^3u_2^2u_3\\
&\hspace{2mm}
 + 
        7776u_1u_2^3u_3 + 19440u_1^4u_3^2 + 26244u_1^2u_2u_3^2 + 
    2916u_2^2u_3^2 + 7290u_1u_3^3 + 9504u_1^6u_4 + 38880u_1^4u_2u_4\\
&\hspace{2mm} - 
        11664u_1^2u_2^2u_4 - 3888u_2^3u_4 + 25272u_1^3u_3u_4 + 
    34992u_1u_2u_3u_4 + 8748u_3^2u_4 + 5832u_1^2u_4^2 + 
    8748u_2u_4^2 \\
&\hspace{2mm}+ 11664u_1^5u_5 + 
        25272u_1^3u_2u_5 + 17496u_1u_2^2u_5 + 26244u_1^2u_3u_5 + 
    8748u_2u_3u_5 + 26244u_1u_4u_5\\
&\hspace{2mm} + 6561u_5^2 + 6561u_3u_6)/6561, \\
h_4&\hspace{-2.5mm}=(1152u_1^{11} + 832u_1^9u_2 + 23616u_1^7u_2^2 - 
    13824u_1^5u_2^3 + 34560u_1^3u_2^4 + 7776u_1^8u_3 + 
    42768u_1^6u_2u_3 + 
        16200u_1^4u_2^2u_3 \\
&\hspace{2mm}+ 7776u_1^2u_2^3u_3 - 5832u_2^4u_3 + 
    17496u_1^5u_3^2 + 27216u_1^3u_2u_3^2 + 23328u_1u_2^2u_3^2 + 
    13122u_1^2u_3^3 + 
        2916u_2u_3^3 \\
&\hspace{2mm}+ 12960u_1^7u_4 - 25920u_1^5u_2u_4 + 
    103680u_1^3u_2^2u_4 - 31104u_1u_2^3u_4 + 19440u_1^4u_3u_4 + 
    69984u_1^2u_2u_3u_4\\
&\hspace{2mm} - 
        5832u_2^2u_3u_4 + 17496u_1u_3^2u_4 + 38880u_1^3u_4^2 - 
    46656u_1u_2u_4^2 - 4374u_3u_4^2 + 2592u_1^6u_5 + 
    21384u_1^4u_2u_5\\
&\hspace{2mm} - 
        5832u_1^2u_2^2u_5 + 11664u_2^3u_5 + 26244u_1^3u_3u_5 + 
    17496u_1u_2u_3u_5 + 6561u_3^2u_5 - 8748u_1^2u_4u_5 + 
    17496u_2u_4u_5 \\
&\hspace{2mm}+ 
        13122u_1u_5^2 + 13122u_4u_6)/13122,\\
h_5&\hspace{-2.5mm}= (10496u_1^{12} + 
    70656u_1^{10}u_2 + 86976u_1^8u_2^2 + 233856u_1^6u_2^3 - 
    25920u_1^4u_2^4 - 20736u_2^6 + 
        71808u_1^9u_3 \\
&\hspace{2mm}
+ 264384u_1^7u_2u_3 + 393984u_1^5u_2^2u_3 + 
    129600u_1^3u_2^3u_3 + 93312u_1u_2^4u_3 + 165888u_1^6u_3^2 + 
    408240u_1^4u_2u_3^2\\
&\hspace{2mm} + 
        221616u_1^2u_2^2u_3^2 + 134136u_1^3u_3^3 + 87480u_1u_2u_3^3 + 
    10935u_3^4 + 22464u_1^8u_4 + 209088u_1^6u_2u_4 + 
    38880u_1^4u_2^2u_4 \\
&\hspace{2mm}+ 
        233280u_1^2u_2^3u_4 + 23328u_2^4u_4 + 241056u_1^5u_3u_4 + 
    272160u_1^3u_2u_3u_4 + 139968u_1u_2^2u_3u_4 + 
    209952u_1^2u_3^2u_4\\
&\hspace{2mm} + 
        87480u_2u_3^2u_4 + 19440u_1^4u_4^2 + 349920u_1^2u_2u_4^2 - 
    69984u_2^2u_4^2 + 139968u_1u_3u_4^2 - 17496u_4^3 + 10368u_1^7u_5\\
&\hspace{2mm} + 
        38880u_1^5u_2u_5 + 93312u_1^3u_2^2u_5 - 23328u_1u_2^3u_5 + 
    134136u_1^4u_3u_5 + 227448u_1^2u_2u_3u_5 + 52488u_2^2u_3u_5 \\
&\hspace{2mm}+ 
        104976u_1u_3^2u_5 + 198288u_1^3u_4u_5 + 104976u_1u_2u_4u_5 + 
    78732u_3u_4u_5 + 78732u_1^2u_5^2 + 26244u_2u_5^2\\
&\hspace{2mm} + 39366u_5u_6)/39366, \\
 h_6&\hspace{-2.5mm}= (109056u_1^{14} + 433664u_1^{12}u_2 + 1983744u_1^{10}u_2^2 - 
    400512u_1^8u_2^3 + 2784768u_1^6u_2^4 - 282240u_1^4u_2^5\\
&\hspace{2mm} + 
        967680u_1^2u_2^6 + 207360u_2^7 + 403200u_1^{11}u_3 + 
    2395008u_1^9u_2u_3 + 3709440u_1^7u_2^2u_3 + 
    6096384u_1^5u_2^3u_3\\
&\hspace{2mm} - 846720u_1^3u_2^4u_3 - 
        725760u_1u_2^5u_3 + 1611792u_1^8u_3^2 + 
    4445280u_1^6u_2u_3^2 + 3492720u_1^4u_2^2u_3^2 + 
    1360800u_1^2u_2^3u_3^2\\
&\hspace{2mm} + 462672u_2^4u_3^2 + 
        1496880u_1^5u_3^3 + 2857680u_1^3u_2u_3^3 + 
    734832u_1u_2^2u_3^3 + 489888u_1^2u_3^4 + 91854u_2u_3^4 \\
&\hspace{2mm}+ 
    177408u_1^{10}u_4 + 80640u_1^8u_2u_4 + 
        3499776u_1^6u_2^2u_4 - 3870720u_1^4u_2^3u_4 + 
    2540160u_1^2u_2^4u_4 - 653184u_2^5u_4 \\
&\hspace{2mm}+ 1439424u_1^7u_3u_4 + 
    7039872u_1^5u_2u_3u_4 + 
        3991680u_1^3u_2^2u_3u_4 + 1741824u_1u_2^3u_3u_4 + 
    2721600u_1^4u_3^2u_4 \\
&\hspace{2mm}+ 1388016u_1^2u_2u_3^2u_4 + 
    244944u_2^2u_3^2u_4 + 
        857304u_1u_3^3u_4 + 1874880u_1^6u_4^2 - 
    3175200u_1^4u_2u_4^2\\
&\hspace{2mm} + 5225472u_1^2u_2^2u_4^2 + 
    925344u_1^3u_3u_4^2 + 2939328u_1u_2u_3u_4^2 + 
        306180u_3^2u_4^2 + 1469664u_1^2u_4^3 - 816480u_2u_4^3 \\
&\hspace{2mm}- 
    8064u_1^9u_5 + 254016u_1^7u_2u_5 + 217728u_1^5u_2^2u_5 + 
    362880u_1^3u_2^3u_5 + 
        544320u_1u_2^4u_5 + 707616u_1^6u_3u_5\\
&\hspace{2mm} + 
    1632960u_1^4u_2u_3u_5 + 2122848u_1^2u_2^2u_3u_5 - 
    326592u_2^3u_3u_5 + 1714608u_1^3u_3^2u_5 + 
        1224720u_1u_2u_3^2u_5\\
&\hspace{2mm} + 183708u_3^3u_5 + 816480u_1^5u_4u_5 + 
    4572288u_1^3u_2u_4u_5 - 1959552u_1u_2^2u_4u_5 + 
    2694384u_1^2u_3u_4u_5\\
&\hspace{2mm} + 
        979776u_2u_3u_4u_5 - 244944u_1u_4^2u_5 + 1102248u_1^4u_5^2 + 
    979776u_1^2u_2u_5^2 + 489888u_2^2u_5^2 + 734832u_1u_3u_5^2\\
&\hspace{2mm} + 
    367416u_4u_5^2 + 
        137781u_6^2)/275562.\\
\end{array}
$$
}

Let $w_j=j/7$ for $j=1,2,3,4,5$ and $w_6=1$.																		
Then each $h_j$ is a weighted homogeneous polynomial of $u_1,\dots,u_6$
if  $w_k$ is the weight of $u_k$ $(k=1,2,\cdots,6)$.
We shall define $6\times 6$ matrices $C$ and $T$ in the following way.
The matrix $C=(C_{ij})$ is defined by $C_{ij}=\partial_{u_i}h_j\>(i,j=1,2,\cdots 6)$.
At this moment, we remark  one of the important properties of $C$.
Namely the matrices $\partial_{ u_j}C\> (j=1,\dots,6)$ are commutative
each other.
For $E=\sum_{j=1}^6w_ju_j\partial_{u_j}$,
we set $T=EC$.
Then $T$ is regarded as a Saito matrix of
the polynomial $F=\det(T)$, where  $F$ is a polynomial of $u_1,u_2,u_3,u_4,u_5,u_6$.
It is expected that $T$ is the Saito matrix of ST34 if we regard  $\{u_1,u_2,u_3,u_4,u_5,u_6\}$
as a set of appropriate basic invariants of ST34.
Actually this is the case as is shown in the theorem below:

\begin{theorem}
We regard $u_1,u_2,\cdots,u_6$ as basic invariants of ST34
which are polynomials of $f_1,f_2,\cdots,f_6$
defined by the following identities:

{\footnotesize
\begin{equation}
\label{equation:uuu-fff}
\left\{
\begin{array}{lll}
f_1 &=& k_1u_1,\\
 f_2 &=& 3k_1^2(484u_1^2 + 125u_2)/56, \\
 f_3 &=& 3k_1^3(3840232u_1^3 + 1647300u_1u_2 + 466785u_3)/3136, \\
  f_4&=& -k_1^4(10u_1^2u_2 - 6u_2^2 - 3u_1u_3 + 3u_4)/175616, \\
   f_5 &=& 27
      k_1^5(8037231129928u_1^5 + 7056011495960u_1^3u_2 + 
               923992672120u_1u_2^2 + 4020065168070u_1^2u_3 \\
&&+ 
        320350835370u_2u_3 + 
               1029935376420u_1u_4 + 170567816805u_5)/9834496, \\
f_6 &=& (41535291386925640512u_1^7 + 
      55402268263969532360u_1^5u_2 + 
            17134764422423982880u_1^3u_2^2\\&& + 
      1037814693245737680u_1u_2^3 + 
            38364019383116279580u_1^4u_3 + 
      13847169906634165920u_1^2u_2u_3 \\&&+ 
            402156952864771860u_2^2u_3 + 1820695718372802075u_1u_3^2 + 
            14361697034479016040u_1^3u_4 \\&&+ 
      2535507816645416760u_1u_2u_4 + 
            322231741664151375u_3u_4 + 5379672850162662240u_1^2u_5\\&& + 
            306822518227402695u_2u_5 + 60643659460340565u_6)/120472576,\\
\end{array}
\right.
\end{equation}
}
where $k_1=(64/27)^{1/7}$.
Then $\det(T)$ is the discriminant of ST34 for the basic invariants   $u_1,u_2,\cdots,u_6$.
\end{theorem}

\textit{Proof}.
We first note that by the identities (\ref{equation:uuu-fff}),
$\{u_1,u_2,\cdots,u_6\}$ is a set of basic invariants of ST34.
Since the action of ST34 on the set of hyperplanes fixed by the pseudo-reflections of ST34 are transitive,
$\det(T)$ is irreducible.
Then to prove the theorem, it is sufficient to show that
$\det(T)=0$ if $x_5=x_6=1$.

To accomplish our purpose, we need some preparation.
We define
symmetric polynomials $r_k$ in $x_1,x_2,x_3,x_4$ by
$$
r_{3j}=x_1^{3j}+x_2^{3j}+x_3^{3j}+x_4^{3j}\>(j=1,2,3,4),\>r_4=x_1x_2x_3x_4,
$$
and
$$
{\tilde p}_{3j}=p_{3j}|_{x_5=x_6=1}\>(j=1,2,3,4,5),\quad {\tilde s}_6=s_6|_{x_5=x_6=1}.
$$
Then
$$
{\tilde p}_{3j}=r_{3j}+2\>(j=1,2,3),\quad {\tilde s}_6=r_4.
$$
and
$$
\begin{array}{lll}
{\tilde p}_{12} &=& \frac{1}{6}(12 + r_3^4 - 6r_3^2r_6 + 3r_6^2 + 8r_3r_9 - 
     24r_4^3), \\
 {\tilde p}_{15} &=&\frac{1}{6}(12 + r_3^5 - 5r_3^3r_6 + 5r_3^2r_9 + 5r_6r_9 - 
     30r_3r_4^3).\\
     \end{array}
     $$
Therefore each matrix entry of
${\tilde T}=T|_{x_5=x_6=1}$ is a polynomial of $r_{3j}\>(j=1,2,3)$ and $r_4$.
Using this expression, we conclude $\det({\tilde T})=0$
by direct computation.

\begin{remark}\label{rem:comment}
Kato, Mano and Sekiguchi \cite{KMS1}, \cite{KMS2} formulated a generalization
of Frobenius manifold structure and  among others 
they introduced the notion of flat coordinates for well-generated complex reflection groups.
As to the group ST34, the second author (J.S.) constructed the Saito matrix for the flat coordinate.
As a consequence, determinant expression of the discriminant of ST34 for the flat coordinate 
was established (cf. \cite{S2}).
In the course of the identification of $u_1,u_2,\cdots,u_6$ with basic invariants of ST34,
we use the Saito matrix constructed by Bessis and Michel.
It is worthwhile to mention  the procedure of this identification.
Noting that the weight of $f_j$ is supposed to be $w_j\ (=\deg f_i/42)$,
 we determine the undetermined  constants $c_{ij}$ and $c_0$  by the conditions
 $$
 \begin{array}{lll}
 f_1&=&c_{11}u_1,\\
 f_2&=&c_{21}u_2+c_{22}u_1^2,\\
 f_3&=&c_{31}u_3+c_{32}u_1u_2+c_{33}u_1^3,\\
 f_4&=&c_{41}u_4+c_{42}u_3u_1+c_{43}u_2^2+c_{44}u_2u_1^2+c_{45}u_1^4,\\
 f_5&=&c_{51}u_5+c_{52}u_4u_1+c_{53}u_3u_2+\cdots,\\
 f_6&=&c_{61}u_6+c_{62}u_5u_2+c_{63}u_5u_1^2+\cdots\\
 \end{array}
 $$
and  
$$F=c_0\det(M_{34}).
$$
Solving these equations, we determine the constants 
 $c_{ij}$ and $c_0$. The answer is given 
 above (\ref{equation:uuu-fff}) which  already appeared in \cite{S2}.
 It is underlined here that the proof of Theorem 2 
 is mainly indebted to (\ref{equation:uuu-fff}) and is independent of the result of Bessis and Michel on the
 explicit form of the Saito matrix $M_{34}$.
 But the determination of
  (\ref{equation:uuu-fff}) depends deeply on the result of Bessis and Michel.
\end{remark}

\section{More about the group ST34}

\subsection{Minimal vectors of $\Lambda^{(3)}$ and the reflections of ST34}
We consider the correspondence between  the totality of minimal vectors of $\Lambda^{(3)}$
and that of hyperplanes fixed by pseudo-reflections of ST34.
The group generated by $P_1,P_2,Q_1,R_1,R_2$ is identified with ST33
(the group numbered as 33 in \cite{ST}).
If $Z_{34}$ denotes the center of ST34, then $Z_{34}$ normalizes ST33, $Z_{34}\cap ST33=\{1\}$ and $|Z_{34}|=6$.
As a consequence $H={\rm ST33}\cdot Z_{34}$ is well-defined as a group and $|H/ST33|=6$.

There is a natural map between coset spaces
$$
ST34/ST33\longrightarrow ST34/H,
$$
in other words,
\[
756 \text{ minimal vectors of }\Lambda^{(3)}\stackrel{6:1}{\longrightarrow}
126 \text{ reflections of ST34}.
\]
We shall construct this map concretely.
For an element $\frac{1}{\sqrt{3}}(\theta, -\omega^b\theta,0,0,0,0)$ up to 
$Z_{34}$,
we assign  $x_1-\omega^b x_2=0$ to this vector.
We take up $\frac{1}{\sqrt{3}}(1,\omega^b,\omega^c,\omega^d,\omega^e,\omega^f), b+c+d+e+f\equiv 0\pmod{3}$.
Without loss of generality, we may assume $0 \leq b \leq c \leq d \leq e \leq f \leq 2$.
Under this condition, we solve $b+c+d+e+f\equiv 0\pmod{3}$.
First we have a trivial solution $(b,c,d,e,f)=(0,0,0,0,0)$ which corresponds to 
$x_1+x_2+x_3+x_4+x_5+x_6=0$.

Case $(b,c,d,e,f)=(0,0,0,0,f),f\not= 0$: This has no solution.

Case $(b,c,d,e,f)=(0,0,0,e,f), e\not=0$: If $e=2$ then $f=2$ and this is not the case. 
Let $e=1$. If $f=2$ then this is the case
corresponding to $x_1+x_2+x_3+x_4+\omega x_5+\omega^2 x_6=0$.
The case $f=1$ does not satisfy the condition.

Case $(b,c,d,e,f)=(0,0,d,e,f),d\not=0$: If $d=2$, then $(e,f)=(2,2)$ satisfy the condition
and this corresponds to
$x_1+x_2+x_3+\omega^2 (x_4+x_5+x_6)=0$.
However, this is the same type as the one we shall see below, that is, the case $(b,c,d,e,f)=(0,0,1,1,1)$.
Let $d=1$. If $e=2$, then $f=2$ and this is not the case. 
Let $e=1$. Then $f=2$ is not the case, but the case $f=1$ works.
This corresponds to $x_1+x_2+x_3+\omega (x_4+x_5+x_6)=0$.

Case $(b,c,d,e,f)=(0,c,d,e,f),c\not=0$: If $c=2$, then $(d,e,f)=(2,2,2)$ and this is not the case.
Let $c=1$. If $d=2$, then $(e,f)=(2,2)$ does not satisfy the condition. 
$(d,e,f)=(1,2,2)$ is the case and 
this corresponds to $x_1+x_2 + \omega (x_3+x_3)+\omega^2 (x_5+x_6)=0$.
The cases $(e,f)=(1,2), (1,1)$ do not satisfy the condition.

Case $(b,c,d,e,f), b\not=0$: If $b=2$ then $(c,d,e,f)=(2,2,2,2)$ does not satisfy the condition.
Let $b=1$. If $c=2$, then $(d,e,f)=(2,2,2)$ satisfies the condition and this corresponds to $x_1+\omega x_2+\omega^2(x_3+x_4+x_5+x_6)=0$. However, this type already appeared in Case 
$(b,c,d,e,f)=(0,0,0,0,1,2)$ as type $x_1+x_2+x_3+x_4+\omega x_5+\omega^2 x_6=0$.
Let $c=1$. If $d=2$, then $(e,f)=(2,2)$ do not satisfy the condition.
Let $d=1$.
The case $(e,f)=(2,2)$ does not satisfy the condition.
The case $(e,f)=(1,2)$ satisfies the condition 
corresponding to $x_1+\omega (x_2+x_3+x_4+x_5)+\omega^2x_6=0$, but this case 
already appeared  
in Case $(b,c,d,e,f)=(0,0,0,1,2)$ as type $x_1+x_2+x_3+x_4+\omega x_5+\omega^2 x_6=0$.
The case $(e,f)=(1,1)$ does not satisfy the condition.

We have thus obtained the following proposition.
\begin{proposition}
There is a natural 6 to 1 correspondence between the minimal vectors of $\Lambda^{(3)}$ up to $Z_{34}$
and the hyperplanes fixed by the pseudo-reflections of ST34.
Typical correspondences are given by
\begin{align*}
(\theta,-\omega^b\theta,0,0,0,0)&\rightarrow x_1-\omega^bx_2=0,\\
(1,1,1,1,1,1)&\rightarrow x_1+x_2+x_3+x_4+x_5+x_6=0,\\
(1,1,1,1,\omega,\omega^2)&\rightarrow 
x_1+x_2+x_3+x_4+\omega x_5+\omega^2 x_6=0,\\
(1,1,1,\omega,\omega,\omega)&\rightarrow 
x_1+x_2+x_3+\omega (x_4+x_5+x_6)=0,\\
(1,1,\omega,\omega,\omega^2,\omega^2)&\rightarrow 
x_1+x_2 + \omega (x_3+x_4)+\omega^2 (x_5+x_6)=0.
\end{align*}
\end{proposition}

\subsection{Basic invariants of ST33 and those of ST34}

We now mention the relationship between the basic invariants of the group ST33 and those of ST34.
For this purpose, we recall the  symmetric polynomials $r_k$ in $x_1,x_2,x_3,x_4$ introduced in 
the proof of Theorem 2.
If we let $x_5=x_6=y$ in $m_j\>(j=1,2,3,4,5,7)$ at \S1,
then we get polynomials in $x_1,x_2,x_3,x_4,y$.
Since we have
$$
p_{3j}|_{x_5=x_6=y}=r_{3j}+2y^{3j}\>(j=1,2,3,4,5,6),
\quad
s_6=r_4y^2,
$$
we can regard $m_j$'s as polynomials in $r_3,r_4,r_6,r_9,y$.
This procedure gives us the invariants of ST33.
First we set
$$
\begin{array}{lll}
J_4 &=& 3r_4 - r_3y + y^4, \\
     J_{10} &=& \frac{1}{2}(r_3^2r_4 - 3r_6r_4 + 2r_3^3y - 
       8r_3r_6y + 6r_9y + 36r_4^2y^2 - 8r_3r_4y^3 - 
       8r_3^2y^4 + 8r_6y^4 + 64r_4y^6)
     \end{array}
$$
and 
$$
J_{6k}=m_k|_{x_5=x_6=y}
\>(k=1,2,3,4,5,7).
$$
Then we get
{\footnotesize
$$
\begin{array}{lll}
   J_6 &=& -5r_3^2 + 6r_6 - 180r_4y^2 - 20r_3y^3 - 8y^6, \\
                   J_{12} &=&\frac{1}{34}(475r_3^4 - 540r_3^2r_6 + 2349r_6^2 - 
       2250r_3r_9 + 81000r_4^3 + 103950r_3^2r_4y^2 - 
       80190r_6r_4y^2 \\
          &&+ 30800r_3^3y^3 - 
              83160r_3r_6y^3 + 49500r_9y^3 + 
       1871100r_4^2y^4 + 415800r_3r_4y^5 + 50820r_3^2y^6 \\
       &&- 
       33264r_6y^6 + 255420r_4y^8 + 6380r_3y^9 + 
              3676y^{12}), \\
   J_{18} &=&\frac{1}{1822}(-86585r_3^6 + 1859625r_3^4r_6 - 
       9692055r_3^2r_6^2 + 3982527r_6^3 - 86580r_3^3r_9 \\
       &&+ 
       12427020r_3r_6r_9 - 
              8402130r_9^2 - 80904420r_3^2r_4^3 + 
       6342300r_6r_4^3 - 167241240r_3^4r_4y^2 \\
       &&+ 
       452096640r_3^2r_6r_4y^2 - 13384440r_6^2r_4y^2 - 
              271985040r_3r_9r_4y^2 - 5634849240r_4^4y^2\\
              && - 
       28007160r_3^5y^3 + 97895520r_3^3r_6y^3 - 
       45110520r_3r_6^2y^3 - 
              54455760r_3^2r_9y^3 + 29082240r_6r_9y^3 \\
              &&- 
       5768197920r_3r_4^3y^3 - 6947020080r_3^2r_4^2y^4 + 
       4465941480r_6r_4^2y^4 - 
              3216213000r_3^3r_4y^5 \\
              &&+ 
       7443235800r_3r_6r_4y^5 - 4300536240r_9r_4y^5 - 
       603515640r_3^4y^6 + 1734248880r_3^2r_6y^6 \\
       &&- 
              18044208r_6^2y^6 - 1099360080r_3r_9y^6 - 
       60989423040r_4^3y^6 - 27788080320r_3r_4^2y^7 \\
       &&- 
       5927021100r_3^2r_4y^8 + 
              4572273420r_6r_4y^8 - 1157156000r_3^3y^9 + 
       3124321200r_3r_6y^9 \\
       &&- 2006741880r_9y^9 - 
       18856197360r_4^2y^{10} - 2278916640r_3r_4y^{11}\\
       && - 
              126792120r_3^2y^{12} + 127646064r_6y^{12 }- 
       200564640r_4y^{14} - 5051040r_3y^{15} - 4823744y^{18}).\\
       \end{array}
       $$
}

Referring to the invariants given by Burkhardt, or by direct calculation,
we can show that $J_4,J_6,J_{10},$ $J_{12},J_{18}$ generate the invariant ring of ST34.
(For example, Burkhardt invariants are given in \cite{NS}.)
Moreover, we can see that 
$J_{24},J_{30},J_{42}$ are polynomials of $J_4,J_6,J_{10},$ $J_{12},J_{18}$ as follows.

{\footnotesize
$$
\begin{array}{lll}
J_{24 }&=&\frac{1}{362524562500}(3749217437791J_{12}^2 - 
       72606621297375000J_{10}^2J_4 
       - 548583360913500J_{12}J_4^3\\
       && + 
              24202207099125000J_4^6 
              + 2954928140625J_{18}J_6 - 
       12101103549562500J_{10}J_4^2J_6 - 
       19382588618832J_{12}J_6^2\\
       && + 
              548583360913500J_4^3J_6^2 
              + 13040967602916J_6^4), \\
       J_{30}&=& \frac{1}{86134928412500}(398092491836609625000J_{10}^3 + 
       1638875315504600J_{12}J_{18} - 
              9023429814963151500J_{10}J_{12}J_4^2\\
              && - 
       129489400924462500J_{18}J_4^3 + 
       796184983673219250000J_{10}J_4^5 - 
              3829737121534827J_{12}^2J_6\\
              && - 
       100339389337262625000J_{10}^2J_4J_6 + 
       25366379163880500J_{12}J_4^3J_6 + 
              33446463112420875000J_4^6J_6 \\
              &&+ 
       1580814161617275J_{18}J_6^2 - 
       7699801741247286000J_{10}J_4^2J_6^2 - 
              10484943427856196J_{12}J_6^3\\
              && + 
       104123021760582000J_4^3J_6^3
        + 11181126000681648J_6^5), \\
       J_{42} &=&\frac{1}{
       7401652548230514062500}(417391800833624595374824312500
        J_{10}^3J_{12} 
        + 1797658480288030291174075J_{12}^2
                J_{18} \\
                &&- 
       15459093584438022965939062500J_{10}^2J_{18}J_4 - 
       9460880818895490828496017750J_{10}J_{12}^2J_4^2 \\
       &&+ 
              22245016831722459672845062500000J_{10}^3J_4^3 - 
       252569017560632216590387500J_{12}J_{18}J_4^3\\
       && + 
              1339003983186291610000803375000J_{10}J_{12}J_4^5 + 
       8770901365894143085607812500J_{18}J_4^6 \\
       &&- 
              44490033663444919345690125000000J_{10}J_4^8 - 
       6443664189964048153587159J_{12}^3J_6 \\
       &&+ 
              504384671966073105703125J_{18}^2J_6
               - 
       2721541301555381788577062500J_{10}^2J_{12}J_4J_6\\
       && - 
              2576515597406337160989843750J_{10}J_{18}J_4^2J_6
               + 
       800906441502367667511497250J_{12}^2J_4^3J_6\\
       && + 
              5561254207930614918211265625000J_{10}^2J_4^4J_6 
              - 
       20983038899424542341187062500J_{12}J_4^6J_6 \\
       &&+ 
              1279530628386314942683324125000J_{10}^3J_6^2 
              - 
       370850790506796659157525J_{12}J_{18}J_6^2 \\
       &&- 
       19995403641453544837088839500J_{10}J_{12}J_4^2J_6^2
        - 
       299396852886277645255706250J_{18}J_4^3J_6^2\\
       && + 
              2054840875253587466115493500000J_{10}J_4^5J_6^2 
              - 
       14237902960750897238629098J_{12}^2J_6^3\\
       && - 
              69323880439964411839202625000J_{10}^2J_4J_6^3 
              + 
       1329325545599266972946080500J_{12}J_4^3J_6^3 \\
       &&+ 
              41380309308849671453485500000J_4^6J_6^3
               + 
       63445090000090696420950J_{18}J_6^4\\
       && + 
              17448714170095736727621576000J_{10}J_4^2J_6^4 + 
       30208710160434535723891548J_{12}J_6^5\\
       &&
        - 
              1578266116654724778611484000J_4^3J_6^5 - 
       11514378808918757251753416J_6^7).\\
       \end{array}
       $$
}

\section{Concluding remarks}

The basic invariants of ST34 obtained by Conway and Sloane \cite{CS}
are difficult to deal with because of the lengthy if we express them as polynomials of
$x_1,\cdots,x_6$.
The aim of the present report is to get their reasonable and explicit forms which are applicable
by the use of basic invariants of $G(3,3,6)$.
As one of applications of the result of this paper, 
 it is possible to show a relationship between the corank one subdiagrams of
 ``Dynkin diagram'' for the group ST34 and subfamilies of 1-parameter deformations
 of a family of deformations of $E_7$-singularity
 which is not versal but is related with  an algebraic Frobenius manifold of type $E_7$ 
 constructed in \cite{S1}.
 This suggests a deep connection between the group ST34 and the real reflection group of type $E_7$.



\begin{thebibliography}{99}
\bibitem{BM} D. Bessis and J. Michel: Explicit presentations for
exceptional braid groups. Experimental Math. (2004), {\bf 13}, 257-266.


\bibitem{Bosmaetal}
W. Bosma, J. Cannon, C. Playoust: The Magma algebra system. I. The user language,
J. Symbolic Comput. 24 (1997), no. 3-4, 235-265. 

\bibitem{CS} J. H. Conway and N. J. A. Sloane: The Coxeter-Todd lattice, the Mitchell group, and related sphere packings.
Math. Proc. Camb. Phil. Soc. (1983), {\bf 93}, 421-440.

\bibitem{KMS1} M. Kato, T. Mano and J. Sekiguchi: Flat structures without potentials.
Rev. Roumaine Math. Pures Appl., {\bf 60} (2015), 4, 481-505.

\bibitem{KMS2} M. Kato, T. Mano and J. Sekiguchi: Flat structure on the space
of isomonodromic deformations,
SIGMA 16(2020), 36 pages.


\bibitem{NS}A. Nagano and H. Shiga: Geometric interpretation of Hermitian modular forms via Burkhardt invariants.
Transformation Groups (2022). https://doi.org/10.1007/s00031-021-09681-w

\bibitem{OT} P. Orlik and H. Terao: {\it Arrangements of Hyperplanes}.
Grundlehren der mathematischen Wissenschaften, 300. Springer-Verlag, Berlin, 1992. 


\bibitem{S2} J. Sekiguchi: Solutions to extended WDVV equations; $ST34,\>E_8$ cases.
Rev. Roumaine Math. Pures Appl., {\bf 64} (2019), 4, 565-583.


\bibitem{S1} J. Sekiguchi: The construction problem of algebraic potentials and reflection groups,
preprint.


\bibitem{ST} G. C. Shephard and A. J. Todd: Finite unitary reflection groups. Canadian J. Math. (1954), {\bf 6}, 274-304.

\bibitem{TE} H. Terao and Y. Enta: Basic derivations for $G_{34}$, Appendix in ``Basic derivations
for complex reflections'' (by P. Orlik),
Contemporary Math. (1989), {\bf 90}, 211-228.

\end{thebibliography}
\end{document}